\input amstex
\input Amstex-document.sty

\def\d{\text{\rm d}}

\def\var{\text{\rm Var}}
\def\ent{\text{\rm Ent}}
\font\teneusm=eusm10 \font\seveneusm=eusm7 \font\fiveeusm=eusm5
\newfam\eusmfam
\textfont\eusmfam=\teneusm \scriptfont\eusmfam=\seveneusm
\scriptscriptfont\eusmfam=\fiveeusm
\def\scr#1{{\fam\eusmfam\relax#1}}

\newsymbol\leqs 1336
\newsymbol\geqs 133E

\pageno 41

\def\shorttitle{Ergodic Convergence Rates of Markov Processes}
\def\shortauthor{Mu-Fa Chen}

\topmatter %
\title\nofrills{\boldHuge Ergodic Convergence Rates of Markov Processes---Eigenvalues, Inequalities and Ergodic Theory}
\endtitle

\author \Large Mu-Fa Chen* \endauthor

\thanks *Department of Mathematics, Beijing Normal University, Beijing 100875,
    China. E-mail: mfchen\@bnu.edu.cn, Home page: http://www.bnu.edu.cn/\~{}chenmf/main\_eng.htm
\endthanks

\abstract\nofrills \centerline{\boldnormal Abstract}

\vskip 4.5mm

{\ninepoint This paper consists of four parts. In the first part,
we explain what eigenvalues we are interested in and show the
difficulties of the study on the first (non-trivial) eigenvalue
through
 examples. In the second part, we present some (dual)
variational formulas and explicit bounds for the first eigenvalue
of Laplacian on Riemannian manifolds or Jacobi matrices (Markov
chains). Here, a probabilistic approach---the coupling methods is
adopted. In the third part, we introduce recent lower bounds of
several basic inequalities; these are based on a generalization of
Cheeger's approach which comes from Riemannian geometry. In the
last part, a diagram of nine different types of ergodicity and a
table of explicit criteria for them are presented. These criteria
are motivated by the weighted Hardy inequality which comes from
Harmonic analysis.

\vskip 4.5mm

\noindent {\bf 2000 Mathematics Subject Classification:} 35P15,
47A75, 49R50, 60J99.

\noindent {\bf Keywords and Phrases:} Eigenvalue, Variational
formula, Inequality, Convergence rate, Ergodic theory, Markov
process.}
\endabstract
\endtopmatter

\document

\baselineskip 4.5mm \parindent 8mm \specialhead \noindent
\boldLARGE  I. Introduction
\endspecialhead

We will start by explaining what eigenvalues we are interested in.

\subhead 1.1\;\; Definition\endsubhead  {\it Consider a
birth-death process with a state space $E\!=\!\{0, 1, 2, \cdots,$
$ n\}$ $(n\leqs\infty)$ and an intensity matrix $Q=(q_{ij})$:
$q_{k,k-1}=a_k>0\,(1\leqs k\leqs n)$, $q_{k,k+1}= b_k>0\,(0\leqs
k\leqs n-1)$, $q_{k,k}=-(a_k+b_k)$, and $q_{ij}=0$ for other $i\ne
j$. }

Since the sum of each row equals 0, we have $Q 1 = 0 = 0\cdot 1$.
This means that the $Q$-matrix has an eigenvalue $0$ with an
eigenvector $1$. Next, consider the finite case of $n<\infty$.
Then, the eigenvalues of $-Q$ are discrete:
$0=\lambda_0<\lambda_1\leqs\cdots \leqs \lambda_{n}$. We are
interested in the first (non-trivial) eigenvalue
$\lambda_1=\lambda_1 -\lambda_0$ (also called spectral gap of
$Q$). In the infinite case ($n=\infty$), $\lambda_1$ can be $0$.
Certainly, one can consider a self-adjoint elliptic operator in
$\Bbb R^d$, the Laplacian $\Delta$ on manifolds, or an
infinite-dimensional operator as in the study of interacting
particle systems.

\subhead 1.2\;\; Difficulties\endsubhead To get a concrete feeling
about the difficulties of this topic, let us first look at the
following examples with a  finite state space. When $E=\{0, 1\}$,
it is trivial that $\lambda_1=a_1+b_0$. The result is nice because
when either $a_1$ or $b_0$ increases, so does $\lambda_1$. When
$E=\{0,1,2\}$, we have  four parameters $b_0, b_1, a_1, a_2$ and
$\lambda_1=2^{-1} \big[a_1+a_2+b_0+b_1
-\sqrt{(a_1-a_2+b_0-b_1)^2+4 a_1 b_1} \big].$ When
$E=\{0,1,2,3\}$, we have six parameters: $b_0, b_1, b_2, a_1, a_2,
a_3$. In this case,  the expression for $\lambda_1$ is too lengthy
to write. The roles of the parameters are inter-related in a
complicated manner. Clearly, it is impossible to compute
$\lambda_1$ explicitly when the size of the matrix is greater than
five.

Next, consider the infinite state space $E=\{0, 1, 2, \cdots\}$.
Denote the eigenfunction of $\lambda_1$  by $g$
 and  the degree of $g$  by $D(g)$
when $g$ is polynomial. Three examples of the perturbation of
$\lambda_1$ and $D(g)$ are listed in Table 1.1.
$$\vbox{\tabskip=0pt \offinterlineskip
\halign to 200pt 
{\strut#& \vrule#\tabskip=0.5em plus1em &\hfil # \hfil& \vrule#&
\hfil# \hfil& \vrule#& \hfil# \hfil& \vrule#& \hfil# \hfil&\vrule#
\tabskip=0pt\cr \noalign{\hrule} & & \omit\hidewidth $b_i (i\geqs
0 )$ \hidewidth & & \omit\hidewidth $a_i (i\geqs 1)$ \hidewidth &
& \omit\hidewidth $\lambda_1$ \hidewidth & & \omit\hidewidth
$D(g)$ \hidewidth &\cr\noalign{\hrule} && $i+c (c>0)$ && $2 i$ &&
$ 1$ && $1$
  &\cr \noalign{\hrule}
&& $i+1$ && $2 i+3$ && $2$ && $2$
    &\cr \noalign{\hrule}
&& $i+1$ && $2 i+\big(4+\sqrt{2}\,\big)$ && $3$ && $3$
    &\cr \noalign{\hrule}
}}$$ \centerline{Table 1.1 Three examples of the perturbation of
$\lambda_1$ and $D(g)$}

\smallskip

\flushpar The first line is the well known linear model for which
$\lambda_1=1$, independent of the constant $c>0$, and $g$ is
linear. Keeping the same birth rate, $b_i=i+1$,  changes the death
rate $a_i$ from $2i$ to $2i+3$ (resp. $2 i+4+\sqrt{2}$), which
leads to the change of $\lambda_1$ from one to two (resp. three).
More surprisingly, the eigenfunction $g$ is changed from linear to
quadratic (resp. triple). For the other values of $a_i$ between $2
i $, $2i+3$ and $2 i+4+\sqrt{2}$,
  $\lambda_1$ is unknown since $g$ is non-polynomial.
As seen from these examples,  the first eigenvalue is very
sensitive. Hence, in general, it is very hard to estimate
$\lambda_1$.

In the next section, we find that this topic is  studied
extensively
 in Riemannian geometry.

\specialhead \noindent \boldLARGE I\!I. New variational formula
for the first eigenvalue
\endspecialhead

\subhead{\bf 2.1\;\; Story of estimating $\lambda_1$ in
geometry}\endsubhead At first,  we recall the study of $
\lambda_1$ in geometry.

Consider Laplacian $ \Delta $ on a compact Riemannian manifold
$(M, g)$, where $g$ is the Riemannian metric. The spectrum of $
\Delta $ is discrete: $\cdots \leqs - \lambda_2\leqs -\lambda_1 <
- \lambda_0=0$ (may be repeated). Estimating these eigenvalues $
\lambda_k$ (especially $ \lambda_1$) is very important in  modern
geometry. As far as we know, five  books, excluding those books on
general spectral theory,  have been devoted to this topic: Chavel
(1984), B\'erard (1986), Schoen and Yau (1988), Li (1993) and Ma
(1993). For a  manifold $M$, denote its dimension,  diameter and
 the lower bound of Ricci curvature by $d$, $D$, and $K$
(Ricci$_M \geqs K g$), respectively.
 We are interested in
estimating $ \lambda_1$ in terms of these three geometric
quantities. It is relatively easy to obtain an upper bound by
applying a test function $f\in C^1(M)$ to the classical
variational formula:
$$\lambda_1=\inf\bigg\{
\int_M \| \nabla f\|^2 \d x : \;  f \in  C^1 (M), \; \int f\d x=0,
\; \int f^2\d x =1 \bigg\}, \tag 2.0$$ where ``$\d x$'' is the
Riemannian volume element. To obtain the lower bound, however, is
much harder. In Table 2.1, we list eight of the strongest lower
bounds
 that have been derived in the past, using various
sophisticated methods.
\def\bbg{$$\matrix \text{P. H. B\'erard, G. Besson}\\
\text{\& S. Gallot (1985)}\endmatrix$$}
\def\zhy{$$\matrix\text{J. Q. Zhong \& }\\
\text{H. C. Yang (1984)}\endmatrix$$}
\def\yja{$$\matrix\text{H. C. Yang (1989) \&}\\
\text{F. Jia (1991)}\endmatrix$$}
$$\vbox{\tabskip=0pt \offinterlineskip
\halign to 12.5cm {\strut#& \vrule#\tabskip=0.3em plus 0.2em
&\hfil # \hfil& \vrule#& \hfil# \hfil&\vrule# \tabskip=0pt\cr
\noalign{\hrule} & & \omit\hidewidth {\bf Author(s)} \hidewidth &
& \omit\hidewidth {\bf Lower bound} \hidewidth
&\cr\noalign{\hrule} && A. Lichnerowicz (1958) &&
$\bold{\dfrac{d}{d-1}\, K, \quad K\geqs 0}. \qquad
\qquad\qquad\qquad\qquad\quad\,\; (2.1)$
  &\cr \noalign{\hrule}
&& $\bbg$ && $\bold{ d\,\bigg\{ \dfrac{\int_0^{\pi/2}\cos^{d-1}t
\d t}{ \int_0^{D/2}\cos^{d-1}t \d t }\bigg\}^{2/d}, \quad
K=d-1>0}. \;\,(2.2)$
  &\cr \noalign{\hrule}
&& P. Li \& S. T. Yau (1980) && $\dfrac{\pi^2}{2\,D^2},\quad
K\geqs 0. \qquad \qquad\qquad \qquad\qquad\qquad\;\;\;(2.3)$
  &\cr \noalign{\hrule}
&& $\zhy$ && $\bold{ \dfrac{\pi^2}{D^2}, \quad K\geqs 0}. \quad
\qquad\qquad\qquad\qquad\qquad\qquad\;\;(2.4)$
  &\cr \noalign{\hrule}
&& P. Li \& S. T. Yau (1980) && $\dfrac{1}{D^2 (d - 1)\exp\big[1 +
 \sqrt{1 + 16\alpha^2}\big]},    \quad K\leqs 0. \quad\;\; (2.5)$
  &\cr \noalign{\hrule}
&& K. R. Cai (1991) && $\bold{ \dfrac{\pi^2}{D^2} + K, \quad
K\leqs 0}. \quad\qquad\qquad  \qquad\qquad\qquad\;\;(2.6)$
  &\cr \noalign{\hrule}
&& $\yja$ && $\bold {\dfrac{\pi^2}{D^2} e^{-\alpha},\quad\text{if
}\;d\geqs 5},
          \quad K\leqs 0. \qquad\qquad\qquad\;\;\;(2.7)$
  &\cr \noalign{\hrule}
&& $\yja$ && $\dfrac{\pi^2}{2\,D^2} e^{-\alpha'}, \quad \text{if
}\; 2\leqs d\leqs 4, \quad K\leqs 0, \qquad\quad\;\; \; \,(2.8)$
  &\cr \noalign{\hrule}
}}$$ \centerline{Table 2.1\quad  Eight lower bounds of
$\lambda_1$}

\flushpar In Table 2.1, the two parameters $\alpha$ and $\alpha'$
are defined as $\alpha=D\sqrt{|K|(d-1)}/2$ and $\alpha' =
D\sqrt{|K|((d-1)\vee 2)}/2$. Among these estimates,  five  ((2.1),
(2.2), (2.4), (2.6) and (2.7))
 are  sharp.  The first two
are sharp for the unit sphere in two  or higher dimensions but
fail for the unit circle; the fourth, the sixth, and the seventh
 are all sharp for the unit circle. As seen from this table,
 the picture is now very complete, due to the efforts
of  many geometers in the past 40 years. Our original starting
point is to learn from the geometers and to study their methods,
especially the recent new developments. In the next section,  we
will show that one can
 go in the opposite direction, i.e., studying the first eigenvalue
 by using probabilistic methods. Exceeding our expectations, we find a
general formula for the lower bound.

\subhead{\bf 2.2\;\; New variational formula}\endsubhead Before
stating our new variational formula, we introduce two notations:
$$
C(r)=\text{cosh}^{d-1}\bigg[ \dfrac{r}{2}\sqrt{
\dfrac{-K}{d-1}}\bigg],
           \; r\in (0, D).\qquad
{\Cal F}=\{f\in C[0,D]: f>0 \text{ on } (0,D)\}.$$ Here, we have
used all the three quantities: the dimension $d$,  the diameter
$D$, and the lower bound $K$ of Ricci curvature.

\smallskip

\flushpar{\bf Theorem 2.1[General formula]\,}(Chen \& Wang
(1997a)).
$$\lambda_1\geqs \sup\limits_{f\in {\Cal F}}\inf\limits_{r\in
(0,D)} \dfrac{4 f(r)}{\int_0^r  C(s)^{-1}\d s\int_s^D C(u)f(u)\d
u}=: \xi_1. \tag 2.9$$

The new variational formula has its essential value in estimating
the lower bound. It is a dual of the classical variational formula
in the sense that ``$\inf$'' in (2.0) is replaced by ``$\sup$'' in
(2.9). The classical formula can be traced to Lord S. J. W.
Rayleigh (1877) and E.\,Fischer\,(1905). Noticing that these two
formulas (2.0) and (2.9) look very different, which explains that
 why such a formula (2.9) has never appeared before. This formula
can produce many  new lower bounds. For instance, the one
corresponding to the trivial function $f \equiv 1$ is  non-trivial
in geometry.  Applying the general formula to the test functions
$\sin(\alpha r)$ and $\cosh^{d-1}(\alpha r)\sin(\beta r)$ with
$\alpha=D\sqrt{|K|(d-1)}/2$ and $\beta={\pi}/{(2D)}$,  we obtain
the following:

\smallskip

\flushpar{\bf Corollary 2.2\,}(Chen\,\&\,Wang (1997a)).
$$\align
& \lambda_1 \geqs \frac{dK}{d-1}\bigg\{1-\cos^d\bigg[\frac{D}{2}
 \sqrt{\frac{K}{d-1}}\bigg]\bigg\}^{-1} , \quad d>1, \quad K\geqs 0, \tag 2.10\\
&\lambda_1\geqs
\frac{\pi^2}{D^2}\sqrt{1-\frac{2D^2K}{\pi^4}}\cosh^{1-d}
\bigg[\frac{D}{2} \sqrt{\frac{-K}{d-1}}\bigg], \quad d>1, \;
K\leqs 0. \text{\hskip -1em}\tag 2.11
\endalign$$

Applying this formula to some very complicated test functions, we
can prove the following result:

\smallskip

\flushpar{\bf Corollary 2.3\,}(Chen,\, Scacciatelli and Yao
(2002)).
$$\lambda_1\geqs {\pi^2}/{D^2}+{K}/{2},\qquad K\in \Bbb R. \tag 2.12$$

The corollaries improve all the estimates (2.1)---(2.8).
Especially, (2.10) improves (2.1) and (2.2), (2.11) improves (2.7)
and (2.8), and (2.12) improves (2.3) and (2.6). Moreover, the
linear approximation in (2.12) is optimal in the sense that the
coefficient 1/2 of $K$ is exact.

A test function is indeed a mimic of the eigenfunction,  so it
should be chosen appropriately  in order to obtain good estimates.
A question  arises naturally: does there exist a single
representative test function such that we can avoid the task of
choosing a different test function each time? The answer is
seemingly negative since we have already seen that the eigenvalue
and the eigenfunction are both very sensitive. Surprisingly, the
answer is affirmative. The representative test function, though
very tricky to find,  has a rather simple form:
$f(r)=\sqrt{\int_0^r C(s)^{-1}\d s}$. This is motivated from the
study of the weighted Hardy inequality, a powerful tool in
harmonic analysis (cf. Muckenhoupt (1972), Opic and Kufner
(1990)).
\smallskip

\flushpar{\bf Corollary 2.4\,}(Chen (2000)).\quad {\it For the
lower bound $\xi_1$ of $\lambda_1$ given in Theorem 2.1, we have}
$$\align
&4 \delta^{-1}\geqs \xi_1\geqs \delta^{-1}, \qquad \text{where}  \tag 2.13\\
&\delta=\sup_{r\in (0, D)}\bigg( \int_0^r C(s)^{-1}\d s\bigg)
\bigg(\int_r^D C(s)\d s\bigg), \qquad
C(s)=\cosh^{d-1}\bigg[\frac{s}{2}\sqrt{\frac{-K}{d-1}}\,\bigg].
\endalign$$

Theorem 2.1 and its  corollaries are also valid  for  manifolds
with a convex boundary endowed with the Neumann boundary
condition. In this case, the estimates (2.1)---(2.8) are
conjectured by the geometers to be correct. However, only the
Lichnerowicz's estimate (2.1) was proven by J. F. Escobar in 1990.
The others in (2.2)---(2.8) and furthermore  in (2.10)---(2.13)
are all new in geometry.

On the one hand, the proof of this theorem is quite
straightforward, based on the coupling introduced by Kendall
(1986) and Cranston (1991). On the other hand, the derivation of
this general formula requires much effort. The key point is to
find a way to mimic the eigenfunctions. For more details,  refer
to Chen (1997).

Applying similar proof techniques to  general Markov processes, we
also obtain  variational formulas for non-compact manifolds,
elliptic operators in $\Bbb R^d$ (Chen and Wang (1997b)), and
Markov chains (Chen (1996)). It is more difficult to derive the
variational formulas for the elliptic operators and Markov chains
due to the presence of infinite parameters in these cases. In
contrast, there are only three parameters ($d$ , $D$, and $K$) in
the geometric case.
 In fact,
formula (2.9) is a particular example of our general formula
(which is complete in dimensional one) for elliptic operators.

To conclude this part, we return to the matrix case introduced at
the beginning of the paper.

\subhead 2.3\;\; Birth-death processes\endsubhead Let
$b_i>0(i\geqs 0)$ and $a_i>0(i\geqs 1)$ be the birth and death
rates, respectively. Define $\mu_0=1$, $\mu_i= b_0\cdots
b_{i-1}/a_1\cdots a_i\, (i\geqs 1)$. Assume that the process is
non-explosive:

\vskip -0.6truecm

$$\tsize\sum_{k=0}^ \infty (b_k \mu_k)^{-1} \sum_{i=0}^k \mu_i = \infty
\qquad \text{and moreover}\qquad \mu=\sum_i \mu_i< \infty.  \tag
2.14
$$

\vskip -0.3truecm

\flushpar The corresponding Dirichlet form is $D(f)\!=\!\sum_i
\pi_i b_i (f_{i+1}-f_i)^2$, $ \scr D (D)\!=\!\{f\!\in\! L^2(\pi):
D(f)< \infty \}$. Here and in what follows, only the dia\-gonal
elements $D(f)$ are written, but the non-diagonal elements can be
computed from the diagonal ones by using the quadrilateral role.
We then have the classical formula $\lambda_1=\big\{D(f):
\pi(f)=0, \pi\big(f^2\big)=1\big\}$. Define $\scr F'=\{f: f_0=0,
\text{there exists $k\!:1\leqs\!k \!\leqs\!\infty$ so that $f_i=
f_{i\wedge k}$}$ $\text{and $f$ is strictly increasing in $[0,
k]$}\}$, $\scr F'' =\{f: f_0=0, \; f \text{ is strictly
increasing}\}$, and $I_i(f)= [\mu_i b_i (f_{i+1}-f_i)]^{-1}
\sum_{j\geqs i+1} \mu_j f_j. $ Let $\bar f= f-\pi (f)$. Then we
have the following results:

\smallskip

\flushpar{\bf Theorem 2.5}\,(Chen (1996, 2000, 2001))\footnote{Due
to the limitation of the space, the most of the author's papers
during 1993--2001 are not listed in References, the readers are
urged to refer to [11].}.\quad {\it Under $(2.14)$, we have
\roster
\item Dual variational formula. $\inf\limits_{f\in \scr F '}\, \sup\limits_{i \geqs 1} I_i (\bar
f)^{-1} = \lambda_1 =\sup\limits_{f\in \scr F''}\,\inf\limits_{i
\geqs 0} I_i (\bar f)^{-1} $.
\item Explicit estimate. $\mu \delta^{-1}\geqs \lambda_1\geqs (4
\delta)^{-1}$, where $ \delta = \sup\limits_{i\geqs 1}
\sum\limits_{j\leqs i-1} (\mu_j b_j)^{-1}\sum\limits_{j\geqs
i}\mu_j$.
\item Approximation procedure. There exist explicit sequences
$\eta_n'$ and $\eta_n''$ such that $ { \eta_n^\prime}^{-1} \geqs
\lambda_1 \geqs {\eta_n^{\prime\prime}}^{-1} \geqs (4 \delta
)^{-1} $.\endroster}

\vskip -0.5truecm

Here the word ``dual'' means that the upper and lower bounds are
interchangeable  if one exchanges ``$\sup$'' and ``$\inf$''. With
slight modifications, this result is also valid for finite
matrices, refer to Chen (1999).

\specialhead \noindent \boldLARGE I\!I\!I. Basic inequalities and
new forms of Cheeger's constants
\endspecialhead

\subhead{\bf 3.1\;\; Basic inequalities}\endsubhead We now go to a
more general setup. Let $(E, \Cal E , \pi )$ be a probability
space satisfying $\{(x, x): x\in E\}\in {\Cal E} \times {\Cal E}$.
Denote by $L^p(\pi)$ the usual real $L^p$-space with norm $\|
\cdot \|_p$. Write $\|\cdot\|=\|\cdot\|_2$.

For a given Dirichlet form $(D, \Cal D(D))$, the classical
variational formula for the first eigenvalue $\lambda_1$ can be
rewritten in the form of  (3.1) below  with an optimal constant
$C= \lambda_1^{-1}$. From this point of view, it is natural to
study other inequalities. Two additional basic inequalities appear
in  (3.2) and (3.3) below.
$$\align
&\text{{\it Poincar\'e inequality}}: \qquad \var(f)\leqs C D(f),
\quad\quad f\in L^2(\pi),\tag 3.1\\
& \text{{\it Logarithmic Sobolev inequality}}\!:\! \int\!\!
f^2\!\log \!\frac{f^2}{\|f\|^2} \d\pi\! \leqs\!  C D(f), \;
f\!\in\!
L^2(\pi), \text{\hskip -2em}\tag 3.2\\
&\text{{\it Nash inequality}}:\;\;\;
  \quad\quad  \var(f)\leqs C D(f)^{1/p} \|f\|_1^{2/q}, \; f\in
  L^2(\pi), \tag 3.3
\endalign$$
where $\var(f)=\pi(f^2)-\pi(f)^2$, $\pi(f)= \int f\d \pi$, $p\in
(1, \infty) $ and $1/p+1/q=1$. The last two inequalities are due
to Gross (1976) and Nath (1958), respectively.

Our main object is a symmetric (not necessarily Dirichlet) form
$(D, \Cal D (D))$ on $L^2(\pi)$, corresponding to an integral
operator (or symmetric kernel) on $(E, \Cal E )$:
$$D(f)\!=\! \frac{1}{2}\!\int_{E\times E}\!J(\d x, \d
y)[f(y)-f(x)]^2,\quad \Cal D (D)=\{f\in L^2(\pi): D(f)< \infty\},
\tag 3.4$$ where $J$ is a non-negative, symmetric measure having
no charge on the diagonal set $\{(x, x): x\in E\}$. A typical
example is the reversible jump process with a $q$-pair $(q(x),$ $
q(x, \d y))$ and a reversible measure $\pi$. Then $J(\d x, \d y)=
\pi(\d x)q(x, \d y)$.

For the remainder of this part, we restrict our discussions to the
symmetric form of (3.4).

\subhead{\bf 3.2\;\; Status of the research}\endsubhead  An
important topic in this research area is to study under what
conditions on the symmetric measure $J$ do the above inequalities
hold. In contrast with the probabilistic method used in Part
(I\!I), here we adopt a ge\-neralization of Cheeger's method
(1970), which comes from Riemannian geometry. Naturally, we define
$\lambda_1:=\inf\{D(f): \pi (f)=0,\, \|f\|=1\}$. For bounded jump
processes, the fundamental known result is the following:
\smallskip

\flushpar{\bf Theorem 3.1}\, (Lawler \& Sokal (1988)). {\it
$\lambda_1\geqs \dfrac{k^2}{2M}$,  where

 \centerline{ $k=\inf\limits_{\pi
(A)\in (0, 1)} \dfrac{\int_A \pi(\d x) q(x, A^c)}{\pi(A)\wedge
\pi(A^c)}$\qquad and \qquad $M=\sup\limits_{x\in E} q(x)$.}}
\smallskip

In the past years, the theorem has been collected into six books:
Chen (1992), Sinclair (1993), Chung (1997), Saloff-Coste (1997),
Colin de Verdi\`ere (1998), Aldous, D. G. \& Fill, J. A. (1994--).
From the titles of the books, one can see a wide range of the
applications. However,  this result fails for the unbounded
operator. Thus, it has been a challenging open problem in the past
ten years to handle the unbounded case.

As for the logarithmic Sobolev inequality, there have been a large
number of publications in the past twenty years for differential
operators. (For a survey, see Bakry (1992) or Gross (1993)).
Still, there are very limited results for integral ope\-rators.
\smallskip

\subhead{\bf 3.3\;\; New results}\endsubhead Since the symmetric
measure can be unbounded, we choose a symmetric, non-negative
function $r(x, y)$ such that \newline $J^{( \alpha )}(\d x, \d y)
   := I_{\{r(x, y)^{\alpha}>0\}}\dfrac{J(\d x, \d y)}{r(x, y)^{ \alpha }}$
$(\alpha >0)$\quad satisfies\quad $\dfrac{J^{(1)}(\d x, E)}{\pi(\d
x)}\leqs 1$, $\pi$-a.s.\newline For convenience, we use the
convention $J^{(0)}=J$. Corresponding to  the three inequalities
above, we introduce  the following new forms of Cheeger's
constants.
$$\vbox{\tabskip=0pt \offinterlineskip
\halign to 12cm {\strut#& \vrule#\tabskip=0.2em plus 0.2em &\hfil
# \hfil& \vrule#& \hfil# \hfil&\vrule# \tabskip=0pt\cr
\noalign{\hrule} & & \omit\hidewidth {\bf Inequality} \hidewidth &
& \omit\hidewidth {\bf Constant $k^{(\alpha)}$} \hidewidth
&\cr\noalign{\hrule} && Poincar\'e && $\inf\limits_{\pi(A)\in (0,
1)} \dfrac{J^{(\alpha)}(A\times A^c)}{\pi (A)\wedge
\pi(A^c)}\qquad\qquad {(\text{Chen \& Wang}(1998))}$
  &\cr \noalign{\hrule}
&& Log. Sobolev && $\lim\limits_{r\to 0}\inf\limits_{\pi(A)\in(0,
r]}
 \dfrac{J^{(\alpha)}(A\times A^c)}
{\pi (A) \sqrt{\log[e +  \pi(A)^{-1}]}}\qquad {(\text{Wang
(2001a)})}$
 &\cr \noalign{\hrule}
&& Log. Sobolev &&$\lim\limits_{\delta\to \infty}
\inf\limits_{\pi(A)>0} \dfrac{J^{(\alpha)}(A\times A^c) +
\delta\pi(A)} {\pi (A) \sqrt{1- \log \pi(A)}}\qquad\;\;
{(\text{Chen (2000)})}$
  &\cr \noalign{\hrule}
&& Nash&& $\inf\limits_{\pi(A)\in (0,
1)}\dfrac{J^{(\alpha)}(A\times A^c)}
 {[\pi (A)\wedge \pi(A^c)]^{(2 q -3)/(2 q -2)}}\qquad\; {(\text{Chen (1999)})}$
  &\cr \noalign{\hrule}
}}$$\centerline{Table 3.1 \quad New forms of Cheeger's constants}
\smallskip

\flushpar Our main result can be easily stated as follows.
\smallskip

\flushpar{\bf  Theorem 3.2}. $k^{(1/2)}>0\Longrightarrow$ {\it the
corresponding inequality holds}.
\smallskip

In other words, we use  $J^{(1/2)}$ and $J^{(1)}$ to handle the
unbounded $J$. The first two kernels come from the use of Schwarz
inequality. This result is proven in four papers quoted in Table
(3.1). In these papers, some estimates which are sharp or
qualitatively sharp for the upper or lower bounds are also
presented.

\specialhead \noindent \boldLARGE I\!V. New picture of ergodic
theory and explicit criteria
\endspecialhead

\subhead{\bf 4.1\;\; Importance of the inequalities}\endsubhead
Let $(P_t)_{t\geqs 0}$ be the semigroup determined by a Dirichlet
form $(D, \Cal D (D))$. Then, various applications of the
inequalities are based on the following results:
\smallskip

\flushpar{\bf  Theorem 4.1}\,(Liggett (1989), Gross (1976) and
Chen (1999)). { \roster
\item Poincar\'e \it inequality $\Longleftrightarrow \|P_t f -\pi (f)\|^2=\var(P_t f)\leqs \var(f)\exp[-2 \lambda_1 t]$.
\item {\it Logarithmic Sobolev inequality} $\Longrightarrow$
  {\it exponential convergence in entropy}: \linebreak $\ent (P_t f)\leqs \ent (f) \exp[- 2\sigma t]$, where
$\ent(f)=\pi(f\log f)$ $ - \pi(f)\log\|f\|_1$.
\item {\it  Nash inequality} $\Longleftrightarrow \var (P_t f)\leqs C
{\|f\|_1}/{t^{1-q}}$.
\endroster}

In the context of diffusions, one can replace ``$\Longrightarrow$"
by ``$\Longleftrightarrow$" in part (2). Therefore, the above
inequalities describe some type of $L^2$-ergodicity for the
semigroup $(P_t)_{t\geqs 0}$. These inequalities have become
powerful tools in the study on infinite-dimensional mathematics
(phase transitions, for instance) and the effectiveness of random
algorithms.
\smallskip

\subhead{\bf 4.2\;\; Three traditional types of
ergodicity}\endsubhead  The following three types of ergo\-dicity
are well known for Markov processes.
$$\align
&\text{\it Ordinary ergodicity}:\qquad \quad
\lim_{t\to  \infty}\|p_t(x, \cdot)-\pi\|_{\var}=0\\
&\text{\it Exponential ergodicity}: \qquad \|p_t(x,
\cdot)-\pi\|_{\var} \leqs C(x) e^{- \alpha t}\;
  \text{for some }\;\;  \alpha >0\\
&\text{\it Strong ergodicity}:\qquad\qquad
  \lim_{t\to  \infty}\sup_x \|p_t(x, \cdot)-\pi\|_{\var}=0\\
&\text{\hskip 12em} \Longleftrightarrow \lim_{t\to \infty}e^{\beta
t}\sup_x \|p_t(x,\cdot)-\pi\|_{\var}=0\; \text{for some } \beta>0
\endalign$$ where $p_t(x, \d y)$ is
the transition function of the Markov process and
$\|\cdot\|_{\var}$ is the total variation norm. They obey the
following implications:
$$\text{Strong
ergodicity}\Longrightarrow\text{Exponential ergodicity}
\Longrightarrow\text{Ordinary ergodicity}.$$ It is natural to ask
the following question. does there exist any relation between the
above inequalities and the three traditional types of ergodicity?
\smallskip

\subhead{\bf 4.3\;\; New picture of ergodic theory}\endsubhead
\smallskip

\flushpar{\bf Theorem 4.2} (Chen (1999), ...). {\it For reversible
Markov processes with densities, we have the diagram shown in
Figure 4.1.
$$\matrix
\format\c\text{\hskip-3em}& \c\text{\hskip-3em} &\c\\
{ }& \text{Nash inequality}&{ }\\
\text{\hskip 4em}\swarrow\!\!\!\!\swarrow&{ }
 &{\text{\hskip -4em}\searrow\!\!\!\!\searrow }\\
\text{Logarithmic Sobolev inequality}&  \text{\hskip 1em}
 & L^1\text{-exponential convergence}\\
\Downarrow &{  } &\Updownarrow\\
\text{Exponential convergence in entropy}
&\text{\hskip 1em} & \text{Strong ergodicity}\\
\Downarrow &{  } &\Downarrow\\
\text{{\hskip -4em}Poincar\'e inequality\;}\text{\hskip
-3em}&\Longleftrightarrow&
    \text{\hskip 0.1em Exponential ergodicity} \\
{ }&\Downarrow&{ }\\
{ }&\text{$L^2$-algebraic ergodicity}&{ }\\
{ }&\Downarrow&{ }\\
{ }& \text{Ordinary ergodicity}&{ }
\endmatrix$$ \centerline{{\rm Figure 4.1 \quad Diagram of nine types of ergodicity}}

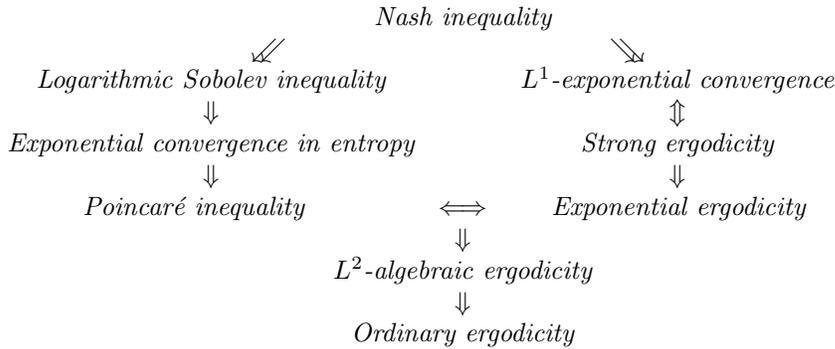
\smallskip

\flushpar In Figure 4.1, $L^2$-{\it algebraic ergodicity} means
that $\var(P_t f)\leqs C V(f) t^{1-q}\,(t>0) $ holds for some $V$
having the properties (cf. Liggett (1991)): $V$ is homogeneous of
degree two (in the sense that $V(cf +d)=c^2 V(f)$ for any
constants $c$ and $d$) and $V(f)<\infty$ for all functions $f$
with finite support. }
\smallskip

The diagram is complete in the following sense: each single-side
implication can not be replaced by double-sides one. Moreover,
strong ergodicity and logarithmic Sobolev inequality (resp.
exponential convergence in entropy) are not comparable. With
exception of the equivalences, all the implications in the diagram
are suitable for more general Markov processes. Clearly, the
diagram extends the ergodic theory of Markov processes.

The diagram was presented in Chen (1999), originally for Markov
chains only. Recently, the equivalence of $L^1$-exponential
convergence and strong ergo\-dicity was mainly proven by Y. H.
Mao. A counter-example of diffusion was constructed by Wang
(2001b) to show that strong ergodicity does not imply exponential
convergence in entropy. For other references and  a detailed proof
of the diagram, refer to Chen (1999).

\subhead{\bf 4.4\;\; Explicit criteria for several types of
ergodicity}\endsubhead As an application of the diagram in Figure
4.1, we obtain a criterion for the exponential ergodicity of
birth-death processes, as listed in Table 4.2. To achieve this, we
use the equivalence of exponential ergodicity and Poincar\'e
inequality, as well as the explicit criterion for Poincar\'e
inequality given in part (3) of Theorem 2.5. This solves a long
standing open problem in the study of Markov chains (cf. Anderson
(1991), \S 6.6 and Chen (1992), \S 4.4).

Next, it is natural to look for some criteria for other types of
ergodicity. To do so, we consider only the one-dimensional case.
Here we focus on  the birth-death processes since the
one-dimensional diffusion processes are in parallel. The criterion
for strong ergodicity was obtained recently by Zhang, Lin and Hou
(2000), and extended by Zhang (2001), using a different approach,
to a larger class of Markov chains. The criteria for logarithmic
Sobolev, Nash inequalities,  and the discrete spectrum (no
continuous spectrum and  all eigenvalues have finite multiplicity)
were obtained by Bobkov and G\"otze (1999) and Mao (2000,
2002a,b), respectively, based on the weighted Hardy inequality
(see also Miclo (1999), Wang (2000), Gong and Wang (2002)). It is
understood now the results can  also be deduced from
generalizations of the variational formulas discussed in this
paper (cf. Chen (2001b)). Finally, we summarize these results in
Theorem 4.3 and Table 4.2. The table is arranged in such an order
that the property in the latter line is stronger than the property
in the former line.  The only exception is that even though the
strong ergodicity is often stronger than the logarithmic Sobolev
inequality, they are not comparable in general, as mentioned in
Part I\!I\!I.

\smallskip

\flushpar{\bf Theorem 4.3}\,(Chen (2001a)). {\it For birth-death
processes with birth rates $ b_i(i\geqs 0)$ and death rates $a_i
(i\geqs 1)$, ten criteria are listed in Table 4.2. Recall the
sequence $(\mu_i)$ defined in Part II and set
$\mu[i,k]=\sum_{i\leqs j\leqs k} \mu_j$. The notion ``$(*) \; \&
\; \cdots$'' appeared in Table 4.2 means that one requires the
uniqueness condition in the first line plus the condition
``$\cdots$''. The notion ``$(\varepsilon)$'' in the last line
means that there is still a small room $(1<q\leqs 2)$ left from
completeness. }

\def\xxp{$$\matrix\text{Exponential ergodicity}\\
\text{$L^2$-exp. convergence}\endmatrix$$}
\def\sxp{$$\matrix\text{Strong ergodicity}\\
\text{$L^1$-exp. convergence}\endmatrix$$} $$\vbox{\tabskip=0pt
\offinterlineskip \halign to 12cm {\strut#& \vrule#\tabskip=0.2em
plus 0.2em &\hfil # \hfil& \vrule#& \hfil# \hfil&\vrule#
\tabskip=0pt\cr \noalign{\hrule} & & \omit\hidewidth {\bf
Property} \hidewidth & & \omit\hidewidth {\bf Criterion}
\hidewidth &\cr\noalign{\hrule} && Uniqueness &&
$\dsize\sum_{n\geqs 0}\frac{1}{\mu_n b_n} \mu[0, n] = \infty \quad
(*)$
  &\cr \noalign{\hrule}
&& Recurrence && $\dsize\sum_{n\geqs 0} \frac{1}{\mu_n
b_n}=\infty$
    &\cr \noalign{\hrule}
&& Ergodicity && $(*) \; \& \;   \mu[0, \infty)< \infty$ &\cr
\noalign{\hrule} && $\xxp$ && $(*) \; \& \; \dsize \sup_{n\geqs 1}
\mu[n, \infty)\!\! \sum_{j\leqs n-1} \frac{1}{\mu_j b_j} <\infty$
    &\cr \noalign{\hrule}
&&\! Discrete spectrum \!\!\! \! && $(*) \; \& \;
\dsize\lim_{n\to\infty}\sup_{k\geqs n\!+\!1} \mu[k, \infty)\!\!
\sum_{n\leqs j\leqs k-1} \frac{1}{\mu_j b_j}=0$
    &\cr \noalign{\hrule}
&& Log. Sobolev inequality && $(*) \; \& \;   \dsize\sup_{n\ge
1}\mu[n,\infty\!) \!\log[
 \mu[n,\infty\!)^{\!-\!1}]\!\!\sum_{j\leqs n\!-\!1}\!
    \! \frac{1}{\mu_j b_j}\!\!<\!\infty $
    &\cr \noalign{\hrule}
&& $\sxp$ && $\!\!(*) \; \& \;   \! \dsize\sum_{n\ge
0}\!\frac{1}{\mu_n b_n}\! \mu[n\!\!+\!\!1,\infty\!)\!= \dsize\!\!
\sum_{n\geqs 1}\!\mu_n \!\!\!\sum_{j\leqs n\!-\!1}\!\!
\frac{1}{\mu_j b_j}\!\!<\!\infty\!\!\!\!\!$
   &\cr \noalign{\hrule}
&& Nash inequality && $(*) \; \& \;   \dsize\sup_{n\geqs 1}
\mu[n,\infty)^{\!(q-2)/(q-1)}\!\!\!
   \sum_{j\leqs n-1} \frac{1}{\mu_j b_j}\!<\!\infty\; (\varepsilon) $
    &\cr \noalign{\hrule}
}}$$ \centerline{Table 4.2 \quad Ten criteria for birth-death
processes} \specialhead \noindent \boldLARGE References
\endspecialhead

\widestnumber\key{100}

\ref \key 1 \by Anderson, W. J. (1991) \book Continuous-Time Markov Chains \publ Springer Series in Statistics
\endref

\ref \key 2 \by Aldous, D. G. \& Fill, J. A. (1994--) \book Reversible Markov Chains and  Random Walks on Graphs,
\rm www. stat.Berkeley.edu/users/aldous/book.html \endref

\ref \key 3 \by Bakry, D. (1992) \paper L'hypercontractivit\'e et
son utilisation en th\'eorie des semigroupes \jour LNM,
Sprin\-ger-Verlag, {\bf 1581}\endref

\ref  \key 4\by B\'erard, P. H. (1986) \book Spectral Geometry:
Direct and Inverse Problem \publ LNM. vol. 1207,
Springer-Verlag\endref

\ref  \key 5\by B\'erard, P. H., Besson, G., Gallot, S.
(1985)\paper Sur une in\'equalit\'e isop\'erim\'etrique qui
g\'en\'eralise celle de Paul L\'evy-Gromov \jour {\it  Invent.
Math.} 80: 295--308\endref

\ref  \key 6\by Bobkov, S. G. and G\"otze, F. (1999)\paper
Exponential integrability and transportation cost related to
       logarithmic Sobolev inequalities
\jour J. Funct. Anal. 163, 1--28
\endref

\ref \key 7 \by Cai, K. R. (1991)\paper Estimate on lower bound of
the first eigenvalue of a compact Riemannian manifold \jour  Chin.
Ann. of Math. 12(B):3, 267--271\endref

\ref \key 8 \by Chavel, I. (1984)\book Eigenvalues in Riemannian
Geometry \publ  Academic Press \endref

\ref  \key 9\by Cheeger, J. (1970)\paper A lower bound for the
smallest eigenvalue of the Laplacian \jour {\it  Problems in
analysis, a symposium in honor of S. Bochner}, 195--199, Princeton
U. Press, Princeton \endref

\widestnumber\key{100} \ref  \key 10\by Chen, M. F. (1992) \book
From Markov Chains to Non-Equilibrium Particle Systems \publ World
Scientific \endref

\ref \key 11\by Chen, M. F. (1993--2001)\book Ergodic Convergence
Rates of Markov Processes {\rm (Collection of papers):
www.bnu.edu.cn/\~{}chenmf/main\_eng.htm}\endref

\ref \key 12 \by Chen, M. F. (2001a) \paper Explicit criteria for
several types of ergodicity \jour Chin. J. Appl. Prob. Stat. 17:2,
1--8\endref

\ref  \key 13\by Chen, M. F. (2001b) \paper Variational formulas
of Poincar\'e-type inequalities in Banach spaces of functions on
the line\jour to appear in Acta Math. Sin. Eng. Ver
\endref

\ref \key 14\by Chen, M. F., Scacciatelli, E. and Yao, Liang
(2002)\paper Linear approximation of the first eigenvalue on
compact manifolds \jour Sci. Sin. (A) 45:4, 450--461  \endref

\ref  \key 15\by Chen, M. F. and Wang, F. Y. (1997a) \paper
General formula for lower bound of the first eigenvalue on
Riemannian manifolds \jour Sci. Sin.
    40:4, 384--394
\endref

\ref  \key 16\by Chen, M. F. and Wang, F. Y. (1997b) \paper
Estimation of spectral gap for elliptic operators \jour Trans.
Amer. Math. Soc. 349, 1239-1267
\endref

\ref  \key 17\by Chen, M. F.  and Wang, F. Y. (1998) \paper
Cheeger's inequalities for general symmetric forms and existence
criteria for spectral gap \jour Abstract. Chin. Sci. Bulletin
43:18, 1516--1519. Ann. Prob. 2000, 28:1, 235--257
\endref

\ref  \key 18\by Chung, F. R. K. (1997) \book Spectral Graph
Theory \publ CBMS, {\bf 92}, AMS, Providence, Rhode Island
\endref

\ref  \key 19\by Colin de Verdi\`ere, Y. (1998) \book Spectres de
Graphes \publ Publ. Soc. Math. France\endref

\ref  \key 20\by Cranston, M. (1991)\paper Gradient estimates on
manifolds using coupling \jour J. Funct. Anal. 99, 110--124\endref

\ref  \key 21\by Escobar, J. F. (1990)\paper Uniqueness theorems
on conformal deformation of me\-trics, Sobolev inequalities, and
an eigenvalue estimate \jour Comm. Pure and Appl. Math. XLI\!I\!I:
857--883\endref

\ref \key 22\by Gong, F. Z. Wang, F. Y.(2002) \paper Functional
inequalities for uniformly integrable semigroups and application
to essential spectrums \jour Forum Math. 14, 293--313
\endref

\ref  \key 23\by Gross, L. (1976) \paper Logarithmic Sobolev
inequalities \jour Amer. J. Math. 97, 1061--1083
\endref

\ref  \key 24\by Gross, L. (1993) \paper Logarithmic Sobolev
inequalities and contractivity properties of semigroups \jour LNM
{\bf 1563}, Sprin\-ger-Verlag\endref

\ref \key 25 \by Jerrum, M. R. and Sinclair, A. J. (1989) \paper
Approximating the permanent \jour SIAM J. Comput.18,
1149--1178\endref

\ref  \key 26\by Jia, F. (1991) \paper Estimate of the first
eigenvalue of a compact Riemannian ma\-nifold with Ricci curvature
bounded below by a negative constant \text{\rm (In Chinese)} \jour
Chin. Ann. Math. 12(A):4, 496--502\endref

\ref \key 27 \by Kendall, W. (1986) \paper Nonnegative Ricci
curvature and the Brownian coupling property \jour Stochastics 19,
111--129\endref

\ref  \key 28\by Lawler, G. F. and Sokal, A. D.(1988) \paper
Bounds on the $L^2$ spectrum for Markov chain and Markov
processes: a generalization of Cheeger's inequality \jour Trans.
Amer. Math. Soc.309, 557--580\endref

\ref  \key 29\by Li, P. (1993) \book Lecture Notes on Geometric
Analysis \publ Seoul National U., Korea\endref

\ref  \key 30\by Li, P. and Yau, S. T. (1980)\paper Estimates of
eigenvalue of a compact Riemannian manifold \jour Ann. Math. Soc.
Proc. Symp. Pure Math. 36, 205--240\endref

\ref  \key 31\by Lichnerowicz, A. (1958)\book Geometrie des
Groupes des Transformationes \publ Dunod \endref

\ref \key 32\by Liggett, T. M. (1989) \paper Exponential $L_2$
convergence of attractive reversible nea\-rest particle systems
\jour Ann. Prob. 17, 403-432
\endref

\ref  \key 33\by Liggett, T. M. (1991) \paper $L_2$ rates of
convergence for attractive reversible nearest particle systems:
the critical case \jour Ann. Prob. 19:3, 935--959
\endref

\ref \key 34 \by Ma, C. Y. (1993) \book The Spectrum of Riemannian
Manifolds \text{\rm (In Chinese)} \publ  Press of Nanjing U.,
Nanjing\endref

\ref  \key 35\by Mao, Y. H. (2000) \paper On empty essential
spectrum for Markov processes in dimension one  \jour preprint
\endref

\ref  \key 36\by Mao, Y. H. (2002a)\paper The logarithmic Sobolev
inequalities for birth-death process and diffusion process on the
line \jour Chin.  J. Appl. Prob. Statis., 18(1), 94--100
\endref

\ref \key 37 \by Mao, Y. H. (2002b)\paper Nash inequalities for
Markov processes in dimension one \jour Acta. Math. Sin. Eng.
Ser., 18(1), 147--156
\endref

\ref  \key 38\by Miclo, L. (1999) \paper An example of application
of discrete Hardy's inequalities \jour Markov Processes Relat.
Fields 5, 319--330
\endref

\ref  \key 39\by Muckenhoupt, B. (1972)\paper Hardy's inequality
with weights \jour Studia Math. XLIV, 31--38
\endref

\ref \key 40\by Nash, J. (1958) \paper Continuity of solutions of
parabolic and elliptic equations \jour Amer. J. Math., 80,
931--954\endref

\ref  \key 41\by Opic, B. and Kufner, A. (1990)  \book Hardy-type
Inequalities \publ Longman, New York
\endref

\ref \key 42 \by Saloff-Coste, L. (1997) \paper Lectures on finite
Markov chains \jour LNM {\bf 1665}, 301--413,
Springer-Verlag\endref

\ref  \key 43\by Schoen, R. and Yau, S. T. (1988) \book
Differential Geometry \text{\rm (In Chinese)} \publ Science Press,
Beijing, China
\endref

\ref  \key 44\by Sinclair, A. (1993) \book Algorithms for Random
Generation and Counting: A Markov Chain Approach \publ
Birkh\"auser\endref

\ref \key 45\by Wang, F. Y. (2000)  \paper Functional inequalities
for empty essential spectrum \jour J. Funct. Anal. 170, 219--245
\endref

\ref  \key 46\by Wang, F. Y. (2001a)\paper Sobolev type
inequalities for general symmetric forms \jour Proc. Amer. Math.
Soc. 128:12, 3675--3682\endref

\ref \key 47\by Wang, F. Y. (2001b)  \paper Convergence rates of
Markov semigroups in probability distances \jour preprint
\endref

\ref  \key 48\by Yang, H. C. (1989)\paper Estimate of the first
eigenvalue of a compact Riemannian manifold with Ricci curvature
bounded below by a negative constant \text{\rm (In Chi-} \text{\rm
ne\-se)} \jour Sci. Sin.(A) 32:7, 698--700\endref

\ref  \key 49\by Zhang, H. J., Lin, X. and Hou, Z. T. (2000)
\paper Uniformly polynomial convergence for standard transition
functions \publ In  ``{\it Birth-death Processes}'' by Hou, Z. T.
et al (2000), Hunan Sci. Press, Hunan
\endref

\ref \key 50 \by Zhang, Y. H. (2001) \paper Strong ergodicity for
continuous-time Markov chains \jour J. Appl. Prob. 38, 270--277
\endref

\ref  \key 51\by Zhong, J. Q. and Yang, H. C. (1984)\paper
Estimates of the first eigenvalue of a compact Riemannian
manifolds \jour Sci. Sin. 27:12, 1251--1265\endref

\enddocument